\documentclass[<options>]{elsarticle}
\usepackage{hyperref}
\usepackage{amsmath,amsthm,amscd,amssymb}
\usepackage{latexsym}
\usepackage{bbm}

\newcommand{\ii}{{\mathrm{i}}}

\newcommand{\lal}{\langle\langle}
\newcommand{\rar}{\rangle\rangle}

\newcommand{\bdone}{{\boldsymbol{1}}}

\renewcommand{\Re}{\operatorname{Re}}
\renewcommand{\Im}{\operatorname{Im}}

\newcommand{\no}{\nonumber}

\newcommand{\tr}{\text{\rm{tr}}}

\newcommand{\beq}{\begin{equation}}
\newcommand{\eeq}{\end{equation}}
\newcommand{\ba}{\begin{align}}
\newcommand{\ea}{\end{align}}

\def \g{{\tt g}}

\def \sur#1#2{\mathrel{\mathop{\kern 0pt#1}\limits^{#2}}}

\newcounter{smalllist}

\allowdisplaybreaks
\numberwithin{equation}{section}

\newtheorem{theorem}{Theorem}[section]
\newtheorem{proposition}[theorem]{Proposition}
\newtheorem{lemma}[theorem]{Lemma}

\theoremstyle{definition}

\theoremstyle{remark}

\title
{A matrix version of a higher-order Szeg{\H o} theorem}
\author
{
Alain Rouault}
\ead{alain.rouault@uvsq.fr}
\address{Laboratoire de Math{\' e}matiques de Versailles, UVSQ, CNRS, Universit\'e Paris-Saclay, 78035-Versailles Cedex France, e-mail: alain.rouault@uvsq.fr}
\date{\today}
\begin{document}



\begin{abstract} 
We extend a higher-order sum rule proved by B. Simon to matrix valued measures on the unit circle and their matrix Verblunsky coefficients. 
\end{abstract}
\begin{keyword}Sum rules, Szeg{\H o}'s theorem, Verblunsky coefficients,  matrix measures on the unit circle,  relative entropy\end{keyword}


\maketitle

\section{Introduction}
A probability measure $\mu$  on the unit circle $\mathbb T$ with infinite support is characterized by its Verblunsky coefficients $(\alpha_j(\mu))_{j \geq 0}$, elemnts in the interior of the unit disc. They are associated  with the Szeg{\H o} recursion of orthogonal polynomials in $L^2(\mathbb T, d\mu)$. A sum rule is an identity between an entropy-like  functional of this measure and a functional of the sequence of its Verblunsky coefficients (for short, we say "V-coefficients" in the sequel). The most famous is  Szeg{\H o}'s theorem.
\begin{theorem}
Let $d\mu= w(\theta) \frac{d\theta}{2\pi} + d\mu_s$ be the Lebesgue decomposition of a probability measure on $\mathbb T$ and let $(\alpha_n)_{n \geq 0}$ its V-coefficients. Then
\begin{align}
\label{SZ}
\int_0^{2\pi} \log w(\theta)  \frac{d\theta}{2\pi} = \sum_0^\infty\log\left(1-|\alpha_k|^2\right)\,. 
\end{align}
where both members can be simultaneously finite or $-\infty$.
\end{theorem}

In his  book \cite{simon05}, B. Simon proved  the following statement (higher-order Szeg{\H o} theorem).

\begin{theorem}[\cite{simon05} Th. 2.8.1]
\label{1.1}
Let $d\mu= w(\theta) \frac{d\theta}{2\pi} + d\mu_s$ be a probability measure on $\mathbb T$ and let $(\alpha_n)_{n \geq 0}$ its V-coefficients. Then
\begin{align}
\notag
\int_0^{2\pi} (1 - \cos \theta) \log w(\theta)  \frac{d\theta}{2\pi} =  &\frac{1}{2}(1 - |1+\alpha_0|^2)  - \frac{1}{2}\sum_0^\infty|\alpha_{k+1} - \alpha_k|^2\\
\label{DH}
&+ \sum_0^\infty \left(\log\left(1-|\alpha_k|^2\right) + |\alpha_k|^2\right)\,,
\end{align}
where both members can be simultaneously finite or $-\infty$.
\end{theorem}
Actually this formula may be written in terms of entropies. For probability measures $\nu$ and $\mu$ on $\mathbb T$, let $\mathcal K (\nu | \mu)$ denote the Kullback-Leibler divergence or relative entropy of $\nu$ with respect to $\mu$: 
\begin{equation}
\label{KL}
{\mathcal K}(\nu\, |\, \mu)= \begin{cases}  \ \displaystyle\int_{\mathbb T}\log\frac{d\nu}{d\mu}\!\ d\nu\;\;& \mbox{if}\ \nu\ \hbox{is absolutely continuous with respect to}\ \mu ,\\
   \  \infty  &  \mbox{otherwise.}
\end{cases}
\end{equation}
Usually, $\mu$ is the reference measure. Here the spectral side will involve the reversed Kullback-Leibler divergence, where $\nu$ is the reference measure and $\mu$ is the argument. In this case, we have that $\mathcal{K}(\nu\,  \, | \, \mu)$ is finite if and only if
\begin{align} \label{eq:KL2}
\int_0^{2\pi} \log w(\theta) \, d\nu(\theta) > - \infty , 
\end{align} 
where $d\mu = w(\theta)d\nu(\theta) + d\mu_s$ is the Lebesgue decomposition of $\mu$ with respect to $\nu$.
If we denote
\label{defl0}
\begin{align}
d\lambda_0(\theta)  = \frac{d\theta}{2\pi} \ , \ d\lambda_1(\theta)  = (1-\cos\theta) \frac{d\theta}{2\pi}
\end{align}
the sum rule (\ref{SZ}) may be written
\begin{align}
\label{SZsr}
\mathcal K\left(\lambda_0\, |\,  \mu\right) = 
 - \sum_0^\infty \log\left(1-|\alpha_k|^2\right)\,,
\end{align}  
and
the sum rule (\ref{DH}) may be written
\begin{align}
\notag
\mathcal K\left(\lambda_1\, |\,  \mu\right) = \mathcal K(\lambda_1 \, |\,  \lambda_0) 
+ \Re\alpha_0 + \frac{|\alpha_0|^2}{2} +\frac{1}{2}\sum_0^\infty|\alpha_{k+1} - \alpha_k|^2\\
\label{KKone}
 - \sum_0^\infty \left(\log\left(1-|\alpha_k|^2\right)+ |\alpha_k|^2\right)
\end{align}  
with 
\[ \mathcal K(\lambda_1 \, |\, \lambda_0) =  \int_0^{2\pi}(1-\cos\theta) \log(1-\cos \theta)\frac{d\theta}{2\pi}  = 1 - \log 2\,.\]
In (\ref{KKone}) both sides may be infinite simultaneously, and they are finite if and only if
\begin{align}
\label{gem0}
\sum_k \alpha_k^4 + |\alpha_{k+1}- \alpha_k|^2 < \infty\,.
\end{align}

Actually, it is easy to include (\ref{KKone}) and (\ref{SZsr}) into a family of sum rules depending on a parameter $\g$ such that $|\g| \leq 1$. Let
\begin{align}\label{KKdefg}
d\lambda_\g(\theta) = (1-\g\cos\theta)\ d\lambda_0(\theta)
\end{align}
(called one single nontrivial moment in \cite{simon05} p. 86). Combining (\ref{KKone}) and Szg{\H o}'s formula, we get, as mentioned in  \cite{GNROPUC} Cor. 5.4 :
\begin{align}\notag
\mathcal K(\lambda_\g\, | \ \mu) &= 
\mathcal K(\lambda_\g\, [ \ \lambda_0)
 + \g \left(\Re \alpha_0 + \frac{|\alpha_0|^2}{2} + \frac{1}{2}\sum_1^\infty |\alpha_k - \alpha_{k-1}|^2\right)\\
\label{KKg}
&+ \sum_0^\infty -\log (1 - |\alpha_k|^2) -\g|\alpha_k|^2\,,
\end{align}
where
\begin{align}
\mathcal K(\lambda_\g\, | \ \lambda_0) = \int (1-\g \cos \theta)\log (1- \g\cos \theta) \frac{d\theta}{2\pi}  = 1 - \sqrt{1 - \g^2} + \log \frac{1 + \sqrt{1-\g^2}}{2}\,.\end{align}
It may be  called GW sum rule, since $\lambda_\g$ is the equilibrium measure in a random matrix model due to Gross and Witten (\cite{GrossW}).

For $\g = 0$, we recover  (\ref{SZ}) formula and when  $\g=1$, we recover (\ref{KKone}).

Simon's proof of Theorem \ref{1.1} (see Sect. 2.8 in \cite{simon05}) was based on the use of the Szeg{\H o} function
\[D(z) = \exp \int_0^{2\pi} \frac{e^{\ii \theta} +z}{e^{\ii \theta} -z}\log w(\theta) \frac{d\theta}{4\pi}\,,\]
the asymptotics of the orthogonal polynomial and Szeg{\H o}'s theorem. 
Later on,  Simon gave another proof of this theorem in Sect.  2.8 of \cite {Simon-newbook}. The new proof  uses a relative Szeg{\H o} function and a step-by-step sum rule provided by the coefficient stripping.  

In a series of papers, Gamboa et al. tackled sum rules on the real line and on the unit circle\footnote{See references in \cite{GNROPUC}.},  on a probabilistic way, using large deviations techniques. The main argument is the uniqueness of the rate function when the large deviations of a random measure are considered under two different encodings. In particular, in \cite{GNROPUC}, they (re)proved  Szeg{\H o}'s theorem as a sum rule, stated a new sum rule for the Hua-Pickrell measure, and 
asked 
for a possible probabilistic proof of the higher-order sum rule quoted above. 
Shortly after, 
 Simon et al. \cite{BSZ} gave that proof.

It turns out that probabilistic tools are robust enough to be extended to matrix measures, which allowed Gamboa et al. to give  a probabilistic proof of the famous matrix  Szeg{\H o}'s theorem of  Delsarte et al.  \cite{delsarte} involving matrix V-coefficients. With the notations of the following section, this theorem says that if $d\mu = w(\theta) d\lambda_0 + d\mu_s$ is a non-trivial matrix-measure, then\footnote{We use $\dagger$ for matrix adjoint, keeping the notation $^*$ for 
reversed polynomials.}
\begin{align}
\int_0^{2\pi}
\log \det w(\theta)  d\lambda_0(\theta) = \sum_0^\infty \log\det (\bdone - \alpha_k\alpha_k^\dagger)\,.
\end{align}

In \cite{GNROPUC} the authors proved also  a matrix version of  the Hua-Pickrell sum rule and 
conjectured a matrix version of the GW sum rule (\ref{KKg}). 

These considerations open the way to two challenges: analytical proof and probabilistic proof.  
The second way seems accessible by combining the machinery of \cite{GNROPUC} and of \cite{BSZ}, i.e. a large deviation for a random measure encoded by its V-coefficients, but it seems more natural to begin with  the first way, which will be done in this note.  
Of course, a possible issue comes from the non-commutativity of the product of matrices, but as usual, the story ends well.

We present the notations and main results  in Sect. \ref{mr}. Theorem \ref{main1} is a matrix-version of (\ref{KKg}) and Prop. \ref{4.1} is a  {\it gem} i.e. a condition of finiteness of the entropy. In Sect. 3, we give the proof of the first result,  involving the coefficient stripping method  and a limiting argument. 
In Sect. 4 we give the proof of the {\it gem}
. Finally Sect. 5 is devoted to the proofs of intermediate results.

\section{Notations and main result}
\subsection{Notations}
\label{mr}
Let us begin with some introductory elements on matrix measures.  For a more detailed exposition, see  \cite{damanik2008analytic} Sect. 1, \cite{derevy} Sect. 4, \cite{GNROPUC} Sect. 6.

Let $p > 1$ be an integer and let $\mathcal M_p$ be the set of complex $p \times p$ matrix measures $\mu$ on $\mathbb T$ which are Hermitian, nonnegative and normalized by $\mu (\mathbb T) = \bdone$ 
 (the $p \times p$ identity matrix). A   
matrix measure is called quasi-scalar if it may be wriiten ${\bf 1}\cdot \sigma$ with $\sigma$ a probability measure on $\mathbb T$. 
A $p \times p$ matrix polynomial is a polynomial with coefficients in $\mathbb C^{p\times p}$.
Given a measure $\mu \in \mathcal M_p$, we define two  inner products on the space of $p \times p$ matrix polynomials by setting
\begin{align*}\langle\langle f, g \rangle\rangle_R = \int f(e^{\ii \theta})^\dagger d\mu(\theta) g(e^{\ii \theta})\\
\langle\langle f, g \rangle\rangle_L = \int g(e^{\ii \theta}) d\mu(\theta) f(e^{\ii \theta})^\dagger
\,.\end{align*}
A sequence of matrix polynomials $(\varphi_j)$ is called right-orthonormal if, and only if, 
\[\lal\varphi_i, \varphi_j\rar_R= \delta_{ij}\bdone \, .\]
A matrix measure is called non-trivial if 
\[\tr\!\ \langle\langle f, f \rangle\rangle_R
 >0\]
for every non-zero polynomial $f$. 
We define the right monic matrix orthogonal polynomials $\Phi_n^R$ by applying the block Gram-Schmidt algorithm  to the sequence $\{\mathbf 1, z\mathbf 1, z^2 \mathbf 1, \dots\}$.  In other words, $\Phi_k^R$ is the unique matrix polynomial 
$\Phi_k^R(z) = z^k \mathbf 1 +$ lower order terms, such that $\lal z^j\mathbf 1, \Phi_k^R\rar_R =0$ for $j=0, \dots, k-1$. 
The normalized orthogonal polynomials are defined by
\[\varphi_0^R = \bdone\ \  ,\ \  \varphi_k^R = \Phi_k^R\kappa_k^R.\]
Here the sequence of $p\times p$ matrices $(\kappa_k^R)$ 
satisfies, for all $k$, the condition $\left(\kappa_k^R\right)^{-1}\kappa_{k+1}^R>0_p$
and is such that the sequence $(\varphi_k^R)$ is orthonormal.
We define the sequence of left-orthonormal polynomials $(\varphi_k^L)$ in the same way except that the above condition is replaced by  $\kappa_{k+1}^L \left(\kappa_k^L\right)^{-1}> 0$. The matrix Szeg\H{o} recursion is then 
\begin{eqnarray}
\label{SzL}
z\varphi_k^L -\rho_k^L\varphi_{k+1}^L &=& \alpha_k^\dagger(\varphi_k^R)^*\\
\label{SzR}
z\varphi_k^R - \varphi_{k+1}^R\rho_k^R &=& (\varphi_k^L)^* \alpha_k^\dagger\,,
\end{eqnarray}
where for all $k\in\mathbb{N}_0$,  
\begin{itemize}
\item $\alpha_k$ belongs to  $\mathbb{B}_p$, the closed unit ball of $\mathbb C^{p\times p}$  defined by
\begin{equation}
\label{defBp}\mathbb{B}_p:=\{M \in \mathbb C^{p\times p}: MM^\dagger \leq \mathbf 1\}\,,\end{equation}
\item $\rho_k^R$ and  $\rho_k^L$ are the so-called defect matrices defined by
\begin{align}
\label{defrho}
\rho^R_k :=  \left(\bdone - \alpha_k\alpha_k^\dagger\right)^{1/2}\  , \  \rho^L_k =    \left(\bdone - \alpha_k^\dagger\alpha_k\right)^{1/2}\,,
\end{align}
\item for a matrix polynomial $P$ with degree $k$, the reversed polynomial $P^*$ is defined by
\[P^*(z) := z^k P(1/\bar z)^\dagger\,.\]
\end{itemize}
Verblunsky's theorem establishes a one-to one correspondance between non-trivial (normalized) matrix measures on $\mathbb T$ and sequences of elements in the interior of $\mathbb B_p$ (Theorem 3.12 in  \cite{damanik2008analytic}).

In an alternative way, these V-coefficients may be introduced as matrix Schur coefficients as follows. Let $F$ be the  Caratheodory (or Herglotz) transform of $\mu$ defined by:
\[F(z) = \int \frac{e^{\ii \theta} + z}{e^{\ii \theta} - z} d\mu(\theta)
 \ , \ z \in \mathbb D = \{ z : |z| < 1\}\,,\]
and $f$ the Schur transform defined by:
\[f(z) = z^{-1}(F(z) -\bdone)(F(z) + \bdone)^{-1}\,,\]
which is equivalent  to
\begin{align}
\label{fF}
F(z) = (\bdone + zf(z))(\bdone-zf(z))^{-1}\,.
\end{align}
The Schur recursion is defined as follows.
At  step $0$ we set
\[\alpha_0 = f(0)\,,\]
which gives the first V-coefficient. We define the defect matrices  (right and left) by
\begin{align}
\label{defect}
\rho_0^R = (\bdone - \alpha_0\alpha_0^\dagger)^{1/2} \ , \ \rho_0^L = (\bdone - \alpha_0^\dagger\alpha_0)^{1/2}\,,
\end{align}
and then, at step $1$ we set
\begin{align}
\label{Schur1}
Sf := f_1 = z^{-1}(\rho_0^R)^{-1}(f(z) - \alpha_0)\left(\bdone- \alpha_0^\dagger f(z)\right)^{-1}\rho_0^L
\end{align}
and the second V-coefficient is
\[\alpha_1 = f_1 (0)\,.\]
The other coefficients are defined with the same algorithm 
\[f_{k+1} = Sf_k, \ \alpha_{k+1} = f_{k+1}(0) , ...\,.\]

The following theorem gives the connection between $F$ and the absolutely continuous part of $\mu$.
\begin{theorem}[\cite{damanik2008analytic} Prop. 3.16]
\label{DPS}
For $z \in \mathbb D$, we have 
\begin{align}
\label{3.52-o}
\Re\!\ F(z) = (\bdone- \bar zf(z)^\dagger)^{-1}(\bdone - |z|^2f(z)^\dagger f(z)) (\bdone -zf(z))^{-1}\,.
\end{align}
and the non-tangential boundary values $\Re F(e^{\ii \theta})$ and $f(e^{\ii \theta})$ exist for a.e.  $\theta$.

If $\mu$ is a normalized matrix measure with Lebesgue decomposition
\[d\mu(\theta) = w(\theta) d\lambda_0(\theta) + d\mu_s (\theta)\]
(where $w$ is a $p\times p$ matrix),  
then for a.e.  $\theta$
\[w(\theta) = \Re F (e^{\ii \theta})\,,\]
and for a.e. $\theta$, $\det w(\theta) \no=0$ if and only if $f(e^{\ii \theta})^\dagger f(e^{\ii \theta}) < \bdone$.
\end{theorem}

\subsection{Main result}
When $\Sigma = \bdone \cdot \sigma$ is a pseudo-scalar measure and $d\mu(\theta) = h(\theta) d\sigma(\theta) + d\mu_s (\theta)$, we define the relative entropy
\begin{align}
\mathcal K(\Sigma\!\ | \ \mu) = - \int_\mathbb T \log\det h(\theta) d\sigma(\theta)\,.
\end{align}
We will consider two reference measures:
\begin{align}
d\Lambda_0 (\theta) = \bdone\cdot  d\lambda_0(\theta) \ , \ d\Lambda_\g (\theta) = \bdone\cdot  d\lambda_\g(\theta)\,.
\end{align}

Our main result is the following.

\begin{theorem} For $|\g| \leq 1$, let $d\mu(\theta) = w(\theta) d\lambda_0(\theta) + d\mu_s(\theta)$ be a non-trivial matrix measure, then
\label{main1}
\begin{align}
\label{mainf1}
\int_0^{2\pi}  (1 - \g \cos \theta) \log\det w(\theta) d\lambda_0(\theta)= \sum_0^\infty \log\det (\bdone - \alpha_k\alpha_k^\dagger) - \g T(\alpha_0, \alpha_1, \cdots)
\end{align}
with
\begin{align}
\label{defT}
T(\alpha_0, \alpha_1, \cdots) := \Re\!\ \tr\!\ (\alpha_0 - \sum_0^\infty \alpha_k\alpha_{k+1}^\dagger)\,, 
\end{align}
or in an equivalent form
\begin{align}
\label{mainf2}
\mathcal K( \Lambda_\g \, | \  \mu) = \mathcal K(\lambda_\g \, | \ \lambda_0) -\sum_0^\infty \log\det (\bdone - \alpha_k\alpha_k^\dagger) + \g  T(\alpha_0, \alpha_1, \cdots)\,.
\end{align}

In (\ref{mainf2}), both sides, which are nonnegative, may be simultaneously infinite.

\end{theorem}

 It is 
exactly  Conjecture 6.11 1. in \cite{GNROPUC}.  For $\g = 0$, we recover of course the matrix Szeg{\H o} formula. 

The right hand side may also be written
\begin{align}
\notag
T(\alpha_0, \alpha_1, \cdots) = \Re \tr\!\ \alpha_0
&+ \frac{1}{2} \tr\!\ \alpha_0\alpha_0^\dagger\\
\label{otherT}
&+ \frac{1}{2}\sum_0^\infty \tr\!\ (\alpha_k - \alpha_{k+1})(\alpha_k^\dagger - \alpha_{k+1}^\dagger) - \sum_0^\infty \tr\!\ \alpha_k\alpha_k^\dagger\,.
\end{align}
According to the definition of B. Simon \cite{Simon-newbook}, the {\it gems} are equivalent conditions for the finiteness of entropies.
Like in Corollary 5.4 in \cite{GNROPUC}, we have the following result.

\begin{proposition}
\label{4.1}
\begin{enumerate}
\item If $|\g| < 1$,
\begin{align}
\label{gema}\mathcal K(\Lambda_\g \!\ | \ \mu) < \infty\ \Longleftrightarrow \sum_k \tr\!\ \alpha_k\alpha_k^\dagger < \infty\end{align}
\item 
\begin{align}
\label{gemb}
\mathcal K(\Lambda_1\!\ | \mu) < \infty &\Longleftrightarrow \sum_k \tr\!\ (\alpha_k\alpha_k^\dagger)^2  + 
\sum_k\tr\!\ (\alpha_{k+1} - \alpha_k)(\alpha_{k+1}^\dagger  - \alpha_k ^\dagger )< \infty\\
\label{gemc}
\mathcal K(\Lambda_{-1}\!\ | \mu) < \infty &\Longleftrightarrow \sum_k \tr\!\ (\alpha_k\alpha_k^\dagger)^2  + 
\sum_k\tr\!\ (\alpha_{k+1} + \alpha_k)(\alpha_{k+1}^\dagger  +\alpha_k ^\dagger )< \infty\,.
\end{align}
\end{enumerate}
\end{proposition}

\section{Proof of Theorem \ref{main1}} 
We need a preliminary remark to reduce the case $\g <0$ to the case $\g >0$.
\begin{lemma}[Simon \cite{simon05} 3.2.6 and \cite{Simon2} 9.5.28]
If $\mu$ is a non-trivial matrix measure and $\tilde\mu$ is defined by
\[d\tilde\mu (\theta) = \begin{cases} d\mu(\pi + \theta)& \hbox{if} \ \theta \in [0, \pi]\\
d\mu(\theta-\pi )& \hbox{if} \ \theta \in [\pi, 2\pi]
\end{cases}
\]
then
\begin{align}
\label{+pi}
\alpha_k(\tilde \mu) = (-1)^{k+1}\alpha_k(\mu) \ , \ (k \geq 0)\,.
\end{align}
\end{lemma}
If $\g = -\gamma$ with $\gamma > 0$, we have,
\[\int (1 - \g\cos\theta) \log\det w(\theta)d\lambda_0(\theta) = \int (1 - \gamma\cos\theta) \log\det \tilde w(\theta) d\lambda_0(\theta)\,,\] 
where $w$ (resp. $\tilde w$) is the a.c. part of $\mu$ (resp. $\tilde \mu$).

If we take for granted the result for $\gamma$, we get
\begin{align}\notag
 \int (1 - \gamma\cos\theta) \log\det \tilde w(\theta) d\lambda_0(\theta) &=\\  \notag =\sum_0^\infty \log\det (\bdone - 
\alpha_k(\tilde\mu)\alpha_k^\dagger)(\tilde\mu))& - \gamma T(\alpha_0(\tilde \mu), \alpha_1(\tilde\mu), \cdots)
\end{align}
but, it is straightforward to see that from (\ref{defT}) and (\ref{+pi}) 
\begin{align}
 T(\alpha_0(\mu), \alpha_1(\mu), \cdots) = - T(\alpha_0(\tilde \mu), \alpha_1(\tilde\mu), \cdots)
\end{align}
so that (\ref{mainf1}) holds true.

From now on, in this section we assume $0 \leq \g\leq 1$.

If $\mu$ is a probability measure on $\mathbb T$ with V-coefficients $(\alpha_j (\mu))_{j \geq 0}$ and if $N$ is some positive integer, we denote by $\mu_N$ the measure whose V-coefficients are shifted: 
\[\alpha_j(\mu_N) = \alpha_{j+N}(\mu) \ , \ j \geq 0\,.\]
When $\mu$ has a density $w$ with respect to $\Lambda_0$, 
we denote by $w_N$ the density of $\mu_N$.

The key point is the following "recursion" theorem, matrix version of Theorem 2.8.2 in \cite{Simon-newbook}, whose proof is postpone to Sect. \ref{proofs}.
\begin{theorem}
\label{recursion}
If $\det w \not=0$ a.e., we have
\begin{align}
\label{2.29}
\int 
\log\det\left(w(\theta)w_1(\theta)^{-1}\right) d\lambda_\g(\theta)
  = \log \det (\bdone - \alpha_0\alpha_0^\dagger) - \g  \Re\!\ \tr\!\ ( \alpha_0 -  \alpha_1  -  \alpha_1\alpha_0\dagger)\,.
\end{align}
\end{theorem}

This implies that $\det w_1 \not= 0$ a.e. and then we may iterate. We get, for $N > 1$
\begin{align}\label{entrel}\int 
\log\det \left(w(\theta)w_N(\theta)^{-1}\right) d\lambda_\g(\theta) 
= G_N(\mu)
\end{align}
where
\begin{align}
\label{itN} 
G_N(\mu) =-\g \Re \tr\!\ (\alpha_N -\alpha_0) + \g \sum_0^{N-1}\Re \tr\!\ \alpha_k
\alpha_{k+1}^\dagger + \sum_0^{N-1}
\log\det (\bdone -\alpha_k \alpha_k^\dagger)\end{align}
In terms of entropy, we have the equivalent form of (\ref{2.29}):
\begin{align}
\label{KKF}
\mathcal K(\Lambda_\g \!\ | \ \mu_N) -  \mathcal K(\Lambda_\g \!\ | \  \mu) = G_N(\mu)\,.
\end{align}
To look for a limit when $N \to \infty$, we need a careful study of 
 $G_N(\mu)$
. We have
\begin{align}
\label{3.6}
G_N(\mu) = -\g \Re \tr\!\ (\alpha_N -\alpha_0) + \frac{\g}{2} \tr\!\  (\alpha_N
\alpha_N^\dagger - \alpha_0\alpha_0^\dagger) - \sum_0^{N-1} A_k\,,
\end{align}
with
\begin{align}
\label{defAk}A_k := -\log\det (1 -\alpha_k \alpha_k^\dagger) - \g\!\  \tr\!\  \alpha_k
\alpha_{k}^\dagger +\frac{\g}{2}\!\  \tr\!\ (\alpha_{k+1} - \alpha_k) (\alpha_{k+1}^\dagger - \alpha_k)^\dagger\,.\end{align}
For $\alpha\alpha^\dagger < 1$, we have 
\begin{align}\notag
- \log\det (1 - \alpha\alpha^\dagger) = \tr\!\ \alpha\alpha^\dagger + \frac{1}{2}\!\ \tr(\alpha\alpha^\dagger)^2 + R(\alpha)\,,
\end{align}
with
\begin{align}
\label{logdet}
R(\alpha) > 0 \  , \ R(\alpha) =o(\tr\!\ (\alpha\alpha^\dagger)^2)\,.\end{align}
This yields
\begin{align}
\label{cruz}
A_k \geq (1-\g) \tr\!\ \alpha_k\alpha_k^\dagger  + \frac{1}{2}\!\  \tr\!\ (\alpha_k\alpha_k^\dagger)^2+\frac{\g}{2}\!\  \tr\!\ (\alpha_{k+1} - \alpha_k) (\alpha_{k+1}^\dagger - \alpha_k^\dagger)\,.
\end{align}
In particular,  $A_k \geq  0$ for every $k$ (remind that we have assumed $\g \geq 0$), which gives 
\[S_N (\mu) := \sum_0^{N-1} A_k \ \uparrow \ 
S_\infty(\mu) = \sum_0^\infty A_k \leq \infty\,,\]
(this argument of monotonicity is like in Simon \cite{Simon-newbook} Prop. 2.8.6.

The identity (\ref{mainf2}) will be the result of two inequalities.

{\bf A)} The first one uses   
  the Bernstein-Szeg{\H o} approximation of 
 $\mu$. 
We know, from Theorem 3.9 in \cite{damanik2008analytic}, for every $\theta$ and every integer $k$, $\varphi_k^R(e^{\ii \theta})$   is invertible and from Theorem 3.11 of the same article that
 the measure
\begin{align}
d\mu^{(N)}(\theta) = \left[\varphi_{N-1}(e^{\ii \theta})\varphi_{N-1}(e^{\ii \theta})^\dagger\right]^{-1} d\lambda_0(\theta)
\end{align}
satisfies
\begin{align}
\alpha_j(\mu^{(N)}) = \begin{cases}\alpha_j (\mu)\ &\mbox{if} \ 0 \leq j \leq N-1\\
0 &\mbox{if} \ j \geq N\,.
\end{cases}
\end{align}
We have $(\mu^{(N)})_N  
= \Lambda_0$. 
We may apply  (\ref{KKF}) with $\mu = \mu^{(N)}$, which 
 gives
\begin{align*}
\mathcal K(\Lambda_\g \!\ | \  \Lambda_0)
 -  \mathcal K(\Lambda_\g \!\ | \ \mu^{(N)}) = G_N(\mu^{(N)}) = \g \Re\!\ \tr\!\ \alpha_0 - \frac{\g}{2}\!\ \tr\!\ \alpha_0\alpha_0^\dagger - S_N(\mu) 
\,.
\end{align*}
Since $\mu^{(N)}$ converges weakly to $\mu$, the lower semicontinuity of 
 $\mathcal K (\Lambda_\g \!\ | \ \cdot)$ gives
\begin{align}
\notag
\mathcal K(\Lambda_\g \!\ | \ \Lambda_0) -  \mathcal K(\Lambda_\g \!\ | \  \mu)  \geq \mathcal K(\Lambda_\g \!\ | \  \Lambda_0)
 -  \liminf_N \mathcal K(\Lambda_\g \!\ | \  \mu^{(N)})\\
\label{sens1}
 \geq 
\g \Re\!\ \tr\!\ \alpha_0 - \frac{\g}{2}\!\ \tr\!\ \alpha_0\alpha_0^\dagger - S_\infty(\mu)
\geq -\infty\,.
\end{align}

{\bf B)} If $\mathcal K(\Lambda_\g \!\ | \ \mu) = \infty$ the inequality
\begin{align}
\label{reverse}
\mathcal K(\Lambda_\g \!\ | \  \Lambda_0) -  \mathcal K(\Lambda_\g \!\ | \  \mu) \leq \g \Re \tr\!\ \alpha_0 - \frac{\g}{2}\!\ \tr\!\ \alpha_0\alpha_0^\dagger - S_\infty(\mu)
\end{align}
is trivial. 

If  $\mathcal K(\Lambda_\g \!\ | \  \mu) < \infty$, then $\det w(\theta) > 0$ a.e. and then from (\ref{KKF}) we have $\det w_N(\theta) > 0$ a.s. too. 
We want to let 
 $N \to \infty$ in (\ref{KKF}) in order to get
(\ref{reverse}). To begin with, 
let us  prove that
\begin{align}
\label{alphazero}
\lim_N \alpha_N(\mu) = 0\,.
\end{align}
From (\ref{KKF}) we deduce
\begin{align*}
G_N(\mu) \leq K(\Lambda_\g \!\ | \  \mu) < \infty\,,
\end{align*}
and then, since
\[-p \leq -\Re \tr\!\ \alpha_N + 
\frac{1}{2} \tr\!\ \alpha_N\alpha_N^\dagger \leq \frac{3p}{2}\]
($p$ is the dimension)
we have $S_\infty (\mu) <\infty$.

Let us split the study into two cases:
\begin{enumerate}
\item
if $0\leq \g <1$,  $S_\infty(\mu) < \infty$ implies
\begin{align}
\label{A}
\sum_k  \tr\!\ \alpha_k\alpha_k^\dagger < \infty
\end{align}
hence (\ref{alphazero}) holds true.
\item if $\g =1$, we have 
\begin{align}
\label{B}
\sum_k  \tr\!\ (\alpha_k\alpha_k^\dagger)^2+  \tr\!\ (\alpha_{k+1} - \alpha_k) (\alpha_{k+1}^\dagger - \alpha_k)^\dagger < \infty
\end{align}
which in particular implies that (\ref{alphazero}) holds true.
\end{enumerate} 

This result has consequences for both sides of (\ref{KKF}). On the one hand, since for every $j$
\[\lim_N \alpha_j(\mu_N) = \lim_N \alpha_{N+j}(\mu) \rightarrow 0\,,\]
the sequence $(\mu_N)$ converges weakly to $\Lambda_0$, so
 using again the semicontinuity, we get
\begin{align}\notag
\mathcal K(\Lambda_\g \!\ | \  \Lambda_0) -  \mathcal K(\Lambda_\g \!\ | \ \mu) \leq \liminf_N \mathcal K(\Lambda_\g \!\ | \  \mu_N)-  \mathcal K(\Lambda_\g \!\ | \  \mu)\,.
\end{align}
On the other hand, from (\ref{3.6}) 
\begin{align}
\label{liminf}
\lim G_N(\mu)
= \g \Re\!\ \tr\!\ \alpha_0 - \frac{\g}{2}\!\ \tr\!\ \alpha_0\alpha_0^\dagger -\g S_\infty(\mu)\,.
\end{align}
and then  (\ref{reverse}) holds true also in this case.

Gathering (\ref{sens1}) and (\ref{reverse})
 ends the proof of (\ref{mainf2}) hence (\ref{mainf1})  when $0\leq \g \leq 1$.
\qedhere
\medskip

\section{Proof of Proposition \ref{4.1}}

We consider only the case $0 \leq \g\leq 1$, since for $-1< \g <0$ the reduction from $\g <0$ to $\gamma > 0$ as in  the beginning of Sect. 3 leads directly to the result.

We already saw in the above section, that when $\mathcal K(\Lambda_g \!\ | \  \mu) < \infty$ and $0\leq \g \leq 1$, the good conditions  are fulfilled.

Conversely, we consider three cases. 

If $0 \leq \g < 1$ and (\ref{A}) is fulfilled, then 
\[- \sum_k \log\det \alpha_k\alpha_k^\dagger < \infty\]
and since 
\[ \tr\!\ (\alpha_{k+1} - \alpha_k) (\alpha_{k+1}^\dagger - \alpha_k)^\dagger \leq 2( \tr\!\ \alpha_k\alpha_k^\dagger+   \tr\!\ \alpha_{k+1}\alpha_{k+1}^\dagger)\,,\]
the expression $T(\alpha_0, \alpha_1, \cdots)$ in (\ref{otherT}) is well defined and finite, so is the left hand side of (\ref{mainf2}) and then $\mathcal K(\Lambda_g \!\ | \  \mu)$ is finite.  

If $\g =1$, condition (\ref{B}), jointly with (\ref{logdet})
entails that
\[\sum_0^\infty -\log\det(
1- \alpha_k\alpha_k^\dagger) - \tr\!\  \alpha_k\alpha_k^\dagger + \frac{1}{2} \tr (\alpha_k - \alpha_{k+1})(\alpha_k^\dagger - \alpha_{k+1}^\dagger) < \infty\]
and then 
 gathering 
 (\ref{mainf2}) and (\ref{otherT}) show that 
 $\mathcal K(\Lambda_g \!\ | \  \mu)$ is finite.

\section{Proofs of intermediate results}
\label{proofs}

\subsection{Proof of Theorem \ref{recursion}}

To compute the LHS of (\ref{2.29}) we need the values of the  Fourier coefficients : \[\int e^{\ii k\theta} \log\det \left(w(\theta) w_1(\theta)^{-1}\right) \frac{d\theta}{2\pi}\ \ \hbox{for} \  k=-1,0, 1\,.\]
The strategy is to approach  $\log\det \left(w(\theta) w_1(\theta)^{-1}\right)$ by a function of $z= re^{\ii \theta}$, sufficiently smooth to apply Cauchy's formula.

In view of Theorem \ref{DPS}, it is natural to approximate $w(\theta)(w_1(\theta))^{-1}$ 
 by $\Re F(z)(\Re F_1 (z))^{-1}$ with $z = re^{\ii \theta}$. We define the auxiliary matrix function:
\begin{align}
\label{2.6.7n-m}
D_0(z) := (\bdone-zf(z))^{-1}(\bdone-zf_1(z))\left(\rho_0^L\right)^{-1}\left(\bdone - f(z) \alpha_0^\dagger\right)\,.
\end{align}

We need the following formula whose proof  is  postponed in Sect. \ref{proofsl}.
\begin{lemma}
\label{FF1}
\begin{align}
\label{2.6.8-nm}
\det \left(\Re F(z)  \left(\Re F_1(z)\right)^{-1}\right) = \det  (D_0(z)D_0(z)^\dagger) \frac{\det (\bdone - |z|^2f(z)^\dagger f(z))}{\det (\bdone- f(z)^\dagger f(z))}\,.
\end{align}
\end{lemma}

From Theorem \ref{DPS}, for a.e. $\theta$ we have
\begin{align*}
\lim_{r\uparrow 1} \det \left(\Re F(re^{\ii \theta})  \left(\Re F_1(re^{\ii \theta}))\right)^{-1}\right) &= \det \left(w(\theta w_1(\theta)^{-1}\right)\\
\lim_{r\uparrow 1} \frac{\det (\bdone - |r|^2f(re^{\ii \theta})^\dagger f(re^{\ii \theta}))}{\det (\bdone- f(re^{\ii \theta})^\dagger f(re^{\ii \theta}))}&=1\,,
\end{align*}
so that,
\begin{align*}
\det(w(\theta)w_1(\theta)^{-1}) = \lim_{r\uparrow 1} \det(D_0(re^{\ii \theta})D_0(re^{\ii \theta})^\dagger)\,,
\end{align*}
and
 the remaining part of the proof  is based  on the study of $\det  (D_0(z)D_0(z)^\dagger)$. Some properties of $D_0$ are collected in the following lemma, whose proof is also in Sect. \ref{proofsl}.

\begin{lemma}
\label{2.6.2}
The function $\det D_0$ is analytic in $\mathbb D$ and non-vanishing. Moreover
\begin{align}
\label{Hp}
h := 2 \log\det D_0 \in 
 H^2(\mathbb D)\,.
\end{align}
\end{lemma}

Since  $h\in H^2(\mathbb D) \subset H^1(\mathbb D)$, we have
\begin{align}
\label{Fh}
\int e^{-\ii \theta}h(e^{\ii \theta})\frac{d\theta}{2\pi} = h'(0)\ , \ 
 \int h(e^{\ii \theta})\frac{d\theta}{2\pi}= h(0)\ , \
\int e^{\ii \theta}h(e^{\ii \theta})\frac{d\theta}{2\pi} =0\,,
\end{align}
and then
\begin{align}
\label{3Fourier}
\int (1-\g \cos \theta) \Re\!\ h(e^{\ii \theta}) d\lambda_0(\theta) =\Re\!\  h(0) - \frac{\g}{2}\Re\!\ h'(0)\,.
\end{align}
Let us compute $h(0)$ and $h'(0)$. As $|z| \to 0$,
\begin{align}
\label{Taylorff1}
\det (\bdone-zf(z)) = 1- z (\tr\!\ \alpha_0) + O(z^2) \ , \ \det (\bdone -zf_1(z)) = 1 - z(\tr\! \ \alpha_1)  + O(z^2)
\end{align}

 Now, formula (\ref{Schur1}) can be inverted into
\begin{align}
\label{Schurinv}
f(z) = (\rho_0^R)^{-1}(\alpha_0 + zf_1(z))\left(\bdone + z \alpha_0^\dagger f_1(z)\right)^{-1} \rho_0^L\,,
\end{align}
which gives the expansion 
\begin{align*}\notag
f(z) 
&=  (\rho_0^R)^{-1}\left(\alpha_0 + z(\bdone-\alpha_0\alpha_0^\dagger) f_1 (z) + O(z^2)\right)\rho_0^L\\
&= \alpha_0 +  z (\rho_0^R
\alpha_1\rho_0^L) +  O(z^2)\,,
\end{align*}
so that
\begin{align*}\notag
\bdone - f(z)\alpha_0^\dagger &= \bdone- 
\alpha_0
\alpha_0^\dagger - z(\rho_0^R
\alpha_1\rho_0^L\alpha_0^\dagger)   +  O(z^2)
= (\rho_0^R)^{2}- z(\rho_0^R\alpha_1\alpha_0^\dagger \rho_0^R)  +  O(z^2)\\
&= \rho_0^R\left(\bdone - z(\alpha_1\alpha_0^\dagger)  +  O(z^2) \right)\rho_0^R
\end{align*}
and 
\begin{align}\notag
\det\left(\bdone - f(z)\alpha_0^\dagger\right) &=\det (\rho_0^R)^2 \det\left(\bdone - z(\alpha_1\alpha_0^\dagger)  +  O(z^2)
 \right)\\
\label{1-falpha}
&= \det (\rho_0^R)^2 \left(\bdone - z\!\ \tr\!\ (\alpha_1\alpha_0^\dagger) +  O(z^2)\right)\,.
\end{align}
Gathering (\ref{2.6.7n-m}), (\ref{Taylorff1})  and (\ref{1-falpha}) and using $\det \rho_0^R =  \det \rho_0^L$ we get
\begin{align*}
\det D_0 (z) =(\det \rho_0^R) 
\left(1 -z\!\ \tr\!\ \left( \alpha_0 -  \alpha_1  - \alpha_1\alpha_0^\dagger \right)+ O(z^2)\right)
\end{align*}
Coming back to the definition of $h$, we get 
\begin{align*}
h(z) = \log\det (\bdone- \alpha_0\alpha_0^\dagger)  
 - 2 z\!\ \tr\!\ \left( \alpha_0 -  \alpha_1  - \alpha_1\alpha_0^\dagger \right)
 + O(z^2)
\end{align*}
and from (\ref{3Fourier})
\begin{align*}
\int(1 - \g \cos\theta) \Re\!\ h(e^{\ii \theta}) d\lambda_0(\theta) = \log\det (\bdone- \alpha_0\alpha_0^\dagger) + \g\Re\!\ \tr\!\ \left( \alpha_0 -  \alpha_1  - \alpha_1\alpha_0^\dagger \right)\,.
\end{align*}
\qedhere

\subsection{Proof of Lemma \ref{FF1}}
\label{proofsl}

To simplify, we omit the variable $z$ if unnecessary. 
Applying (\ref{3.52-o}) to $F_1$
\begin{align}
\label{3.52-1}
\Re F_1  = (\bdone- \bar zf_1^\dagger)^{-1}(\bdone - |z|^2f_1^\dagger f_1) (1 -zf)^{-1}
\end{align}
  so we need an expression of $\bdone - |z|^2 f_1f_1^\dagger$ as a function of $f$. From (\ref{Schur1}) we get
\begin{align*} |z|^2f_1(z)^\dagger f_1(z) =&\\  \rho_0^L\left(\bdone- f(z)^\dagger \alpha_0 \right)^{-1}&(f(z)^\dagger - \alpha_0^\dagger)(\rho_0^R)^{-2}(f(z) - \alpha_0)\left(\bdone- \alpha_0^\dagger f(z)\right)^{-1}\rho_0^L
\end{align*}
which,
with the help of the trivial identity
\begin{align}
\notag
\bdone &=   \rho_0^L\left(\bdone- f(z)^\dagger \alpha_0 \right)^{-1}\left(\bdone- f(z)^\dagger \alpha_0 \right)(\rho_0^L)^{-2}\left(\bdone- \alpha_0^\dagger f(z)\right)\left(\bdone- \alpha_0^\dagger f(z)\right)^{-1} \rho_0^L \,,
\end{align}
yields
\begin{align}
\notag
\left(\bdone- f^\dagger \alpha_0 \right)(\rho_0^L)^{-1}
\left(\bdone- |z|^2f_1^\dagger f_1\right)(\rho_0^L)^{-1}\left(\bdone- \alpha_0^\dagger f\right)  &= 
\\
\label{5.10} \left(\bdone- f^\dagger \alpha_0 \right)(\rho_0^L)^{-2}\left(\bdone- \alpha_0^\dagger f\right)-
(f^\dagger - \alpha_0^\dagger)(\rho_0^R)^{-2}&(f - \alpha_0)\,.
\end{align}
Now, we use (\ref{defect}) 
and
\[
 (\rho_0^R)^{-2} = \sum_{n\geq 0} (\alpha_0\alpha_0^\dagger)^n \ , \  (\rho_0^L)^{-2}
 = \sum_{n\geq 0} (\alpha_0^\dagger\alpha_0)^n\]
($\alpha_0$ is a contraction). 
Expanding the RHS of (\ref{5.10}) and cancelling terms gives
\[\left(\bdone- f^\dagger \alpha_0 \right)(\rho_0^L)^{-2}\left(\bdone- \alpha_0^\dagger f\right)-
(f^\dagger - \alpha_0^\dagger)(\rho_0^R)^{-2}(f - \alpha_0) = \bdone- f^\dagger f\]
so that
\begin{align}
\label{2.6.6}
\bdone -| z|^2f_1^\dagger f_1 = \rho_0^L\left(\bdone- f^\dagger \alpha_0 \right)^{-1}\left(\bdone- f^\dagger f\right)\left(\bdone- f \alpha_0^\dagger \right)^{-1}\rho_0^L\,.
\end{align}
Plugging into (\ref{3.52-1}) yields
\begin{align*} 
&\Re F_1\\ &= (\bdone - \bar zf_1^\dagger)^{-1} \rho_0^L\left(\bdone- f(z)^\dagger \alpha_0 \right)^{-1}\left(\bdone- f(z)^\dagger f(z)\right)\left(\bdone- f(z) \alpha_0^\dagger \right)^{-1}\rho_0^L(\bdone - zf_1)^{-1}
\end{align*}
and
\begin{align}\nonumber
&(\Re F)  \left(\Re F_1\right)^{-1}= (\bdone- \bar zf^\dagger)^{-1}(\bdone - |z|^2f^\dagger f) (\bdone -zf)^{-1}  \\
\notag &\times(\bdone - zf_1)\left(\rho_0^L\right)^{-1}\left(\bdone- f \alpha_0^\dagger \right)\left(\bdone- f^\dagger f\right)^{-1}
\left(\bdone- f^\dagger \alpha_0 \right)\left(\rho_0^L\right)^{-1}(\bdone - \bar z f_1^\dagger)\\
\notag
&= (\bdone- \bar zf^\dagger)^{-1} (\bdone - |z|^2f^\dagger f) D_0\left(\bdone- f^\dagger f\right)^{-1}D_0^\dagger(1-\bar zf^\dagger)\,.
\end{align}
Then, taking determinants
\begin{align}
\notag
\det \left((\Re F)  \left(\Re F_1\right)^{-1}\right) = \det  (D_0D_0^\dagger) \ \frac{\det (\bdone - |z|^2f^\dagger f)}{\det (\bdone- f^\dagger f)}\,,
\end{align}
ends the proof.

\subsection{Proof of Lemma \ref{2.6.2}} We  repeat here the argument of 
Theorem 2.6.2 in \cite{Simon-newbook} for the sake of completeness. 
For  $z \in \mathbb D$ we have $f(z)f^\dagger(z) < \bdone$ , hence analyticity and non-vanishing are straightforward. Moreover, since  $|\zeta| < 1$ implies $|\arg(1-\zeta)| < \pi/2$,  
we conclude from (\ref{2.6.7n-m}) and (\ref{Hp}) that 
\[ |\Im h|< 3\pi/2\,.\]
Since $|h|^2 - 2(\Im h)^2$ is harmonic we have
\begin{align}\notag\int |h|^2  d\lambda_0- 2 \int (\Im h)^2 d\lambda_0 = |h(0)|^2 - 2 \left(\Im h(0)\right)^2\,, 
\end{align}
and since $h(0) = \log\det(1 - \alpha_0\alpha_0^\dagger)  < 0$, we get
\[\int|h(re^{\ii \theta})|^2 d\lambda_0(\theta) \leq \frac{9\pi^2}{2}+\left(  \log\det(1 - \alpha_0\alpha_0^\dagger)\right)^2\,,\]
which yields
\[\sup_{r <1} \int|h(re^{\ii \theta})|^2 d\lambda_0(\theta)  < \infty\,.\]
\bibliographystyle{plain}
\bibliography{bibil}
\end{document}